\documentclass{amsart}
\usepackage{amssymb}
\usepackage{graphicx}

\newtheorem{thm}{Theorem}
\newtheorem{cor}{Corollary}
\newtheorem{lem}{Lemma}

\newtheorem{rem}{Remark}

\newcommand{\A}{{\mathcal A}}

\newcommand{\ID}{{\mathbb D}}

\newcommand{\D}{{\mathbb D}}

\def\be{\begin{equation}}
\def\ee{\end{equation}}

\newcommand{\bee}{\begin{enumerate}}
\newcommand{\eee}{\end{enumerate}}

\newcommand{\blem}{\begin{lem}}
\newcommand{\elem}{\end{lem}}
\newcommand{\bthm}{\begin{thm}}
\newcommand{\ethm}{\end{thm}}
\newcommand{\bcor}{\begin{cor}}
\newcommand{\ecor}{\end{cor}}
\newcommand{\beg}{\begin{example}}
\newcommand{\eeg}{\end{example}}
\newcommand{\begs}{\begin{examples}}
\newcommand{\eegs}{\end{examples}}
\newcommand{\bdefe}{\begin{defin}}
\newcommand{\edefe}{\end{defin}}
\newcommand{\bprob}{\begin{prob}}
\newcommand{\eprob}{\end{prob}}
\newcommand{\bei}{\begin{itemize}}
\newcommand{\eei}{\end{itemize}}

\newcommand{\bcon}{\begin{conj}}
\newcommand{\econ}{\end{conj}}
\newcommand{\bcons}{\begin{conjs}}
\newcommand{\econs}{\end{conjs}}
\newcommand{\bprop}{\begin{propo}}
\newcommand{\eprop}{\end{propo}}
\newcommand{\br}{\begin{rem}}
\newcommand{\er}{\end{rem}}
\newcommand{\brs}{\begin{rems}}
\newcommand{\ers}{\end{rems}}
\newcommand{\bo}{\begin{obser}}
\newcommand{\eo}{\end{obser}}
\newcommand{\bos}{\begin{obsers}}
\newcommand{\eos}{\end{obsers}}
\newcommand{\bpf}{\begin{pf}}
\newcommand{\epf}{\end{pf}}
\newcommand{\ba}{\begin{array}}
\newcommand{\ea}{\end{array}}
\newcommand{\beq}{\begin{eqnarray}}
\newcommand{\beqq}{\begin{eqnarray*}}
\newcommand{\eeq}{\end{eqnarray}}
\newcommand{\eeqq}{\end{eqnarray*}}

\begin{document}
\bibliographystyle{amsplain}

\title[Improvement of estimates of  modulus of Hankel determinants ]{An improvement of the estimates of the modulus of the Hankel determinants of second and third order for  the class $\boldsymbol{\mathcal{S}}$ of univalent functions}

\author[M. Obradovi\'{c}]{Milutin Obradovi\'{c}}
\address{Department of Mathematics,
Faculty of Civil Engineering, University of Belgrade,
Bulevar Kralja Aleksandra 73, 11000, Belgrade, Serbia}
\email{obrad@grf.bg.ac.rs}

\author[N. Tuneski]{Nikola Tuneski}
\address{Department of Mathematics and Informatics, Faculty of Mechanical Engineering, Ss. Cyril and Methodius
University in Skopje, Karpo\v{s} II b.b., 1000 Skopje, Republic of North Macedonia.}
\email{nikola.tuneski@mf.edu.mk}


\begin{abstract}
Using some properties of the Grunsky coefficients we improve earlier results for upper bounds of the Hankel determinants of the second and third order for the class  $\mathcal{S}$ of univalent functions.
\end{abstract}

\maketitle

\section{Introduction and preliminaries}

Let $\mathcal{A}$ be the class of functions $f$ analytic  in the open unit disc $\D=\{z:|z|<1\}$ and normalised such that $f(0)=f'(0)-1=0$, i.e., of the form $f(z)=z+a_2z^2+a_3z^3+\cdots$, and let its subclass $\mathcal{S}$ consist of univalent functions in the unit disc $\D$. Further, let $\mathcal{S}^{\star}$ and $\mathcal{K}$ denote the subclasses of
${\mathcal A}$ which are starlike and convex in $\ID$, respectively, and
let $\mathcal{U} $ denote the set of all  $f\in {\mathcal A}$ satisfying
$$\left |\left (\frac{z}{f(z)} \right )^{2}f'(z)-1\right | < 1 \quad\quad (z\in \ID).$$
(see  \cite{OP-01,OP-2011}).

\medskip

In the study of the class of univalent functions and its subclasses, a significant topic is finding upper estimates (preferably sharp) of the Hankel determinant, especially of the second and third order, for a function $f$ from $\A$ is defined by 
\be\label{eq 1}
H_{2}(2)= 
\left |
        \begin{array}{cc}
        a_{2} & a_{3}\\
        a_{3}&a_{4}\\
        \end{array}
        \right |
= a_2a_4-a_{3}^2
\ee
and 
\be\label{eq 2}
H_3(1) =  \left |
        \begin{array}{ccc}
        1 & a_2& a_3\\
        a_2 & a_3& a_4\\
        a_3 & a_4& a_5\\
        \end{array}
        \right | = a_3(a_2a_4-a_{3}^2)-a_4(a_4-a_2a_3)+a_5(a_3-a_2^2),
\ee
respectively.

\medskip

Hankel determinants are used for studies in the theory of singularities (see \cite{dienes}), as well, as in the study of power series with integral coefficients. The upper bound of their modulus is of special interest in the theory of univalent functions and for some subclasses of the class $\mathcal{S}$ the sharp estimation of $|H_{2}(2)|$ are known. For example,
for the classes $\mathcal{S}^{\star}$ and $\mathcal{U} $ we have that $|H_{2}(2)|\leq 1$
(see \cite{JHD}, \cite{OT-2019}), while $|H_{2}(2)|\leq \frac{1}{8} $ for the class $\mathcal{K}$ (\cite{JHD}).
The sharp estimate of $H_3(1)$ seams to be more challenging problem and quite few are known.
A review on this can be found in \cite{shi}, while new non-sharp upper bounds for different classes and conjectures about the sharp ones are given in \cite{OT-2019-1}.

\medskip

 In their paper \cite{OT-2021} the authors gave the next upper bound of $|H_{2}(2)|$ and $|H_{3}(1)|$
for the class $\mathcal{S}$:

\bthm\label{19-th 1} For the class $\mathcal{S}$ we have
\[
|H_{2}(2)|\leq A,  \quad\mbox{where}\quad 1\leq A\leq \frac{11}{3}=3,666\ldots
\]
and
\[
|H_{3}(1)|\leq  B, \quad\mbox{where}\quad \frac49\leq B\leq \frac{32+\sqrt{285}}{15} = 3.258796\cdots
\]
\ethm

\medskip

 In this paper we improve these results by proving:
\bthm\label{19-th 2}
For the class $\mathcal{S}$ we have the next estimations:
\begin{itemize}
  \item[$(i)$]  $|H_{2}(2)|\leq 1.3614\ldots ; $
  \item[$(i)$]  $|H_{3}(1)|\leq 1.6787\ldots.  $
\end{itemize}
\ethm

\medskip

The proof of this theorem will make use mainly the notations and results given in the book of N.A. Lebedev (\cite{Lebedev}).

\medskip

For an univalent function $f$ from $\mathcal{S}$ we have
\[
\log\frac{f(t)-f(z)}{t-z}=\sum_{p,q=0}^{\infty}\omega_{p,q}t^{p}z^{q},
\]
where $\omega_{p,q}$ are the so-called Grunsky's coefficients such that $\omega_{p,q}=\omega_{q,p}$. This coefficients satisfy the Grunsky's inequality (\cite{duren,Lebedev}):
\be\label{eq 3}
\sum_{q=1}^{\infty}q \left|\sum_{p=1}^{\infty}\omega_{p,q}x_{p}\right|^{2}\leq \sum_{p=1}^{\infty}\frac{|x_{p}|^{2}}{p},
\ee
where $x_{p}$ are arbitrary complex numbers such that last series converges.

\medskip

Next, it is well-known that if
\be\label{eq 4}
f(z)=z+a_{2}z^{2}+a_{3}z^{3}+\ldots
\ee
belongs to $\mathcal{S}$, then also does 
\be\label{eq 5}
f_{2}(z)=\sqrt{f(z^{2})}=z +c_{3}+c_{5}z^{5}+\ldots.
\ee
Then, for the function $f_{2}$ the appropriate Grunsky's
coefficients are of the form $\omega_{2p-1,2q-1}^{(2)}$ and the inequality \eqref{eq 3} appears to be
\be\label{eq 6}
\sum_{q=1}^{\infty}(2q-1) \left|\sum_{p=1}^{\infty}\omega_{2p-1,2q-1}^{(2)}x_{2p-1}\right|^{2}\leq \sum_{p=1}^{\infty}\frac{|x_{2p-1}|^{2}}{2p-1}.
\ee
Finally, from \cite[p.57]{Lebedev} we have that the coefficients $a_{2}, a_{3}, a_{4}$ of $f$ can be expressed by Grunsky's coefficients  $\omega_{2p-1,2q-1}^{(2)}$ of $f_{2}$ given by \eqref{eq 5} as:
\be\label{eq 7}
\begin{split}
a_{2}&=2\omega _{11},\\
a_{3}&=2\omega_{13}+3\omega_{11}^{2}, \\
a_{4}&=2\omega_{33}+8\omega_{11}\omega_{13}+\frac{10}{3}\omega_{11}^{3}\\
a_{5}&=2\omega_{35}+8\omega_{11}\omega_{33}+5\omega_{13}^{2}+18\omega_{11}^{2}\omega_{13}+\frac{7}{3}\omega_{11}^{4}\\
0&=3\omega_{15}-3\omega_{11}\omega_{13}+\omega_{11}^{3}-3\omega_{33}\\
0&=\omega_{17}-\omega_{35}-\omega_{11}\omega_{33}-\omega_{13}^{2}+\frac{1}{3}\omega_{11}^{4}.
\end{split}
\ee
Here and in the rest of the paper, for simplicity of    the expressions,  we omit upper index "(2)" in $\omega_{2p-1,2q-1}^{(2)}$.

\medskip

We note that in the book \cite{Lebedev} there exists a typing mistake for the coefficient $a_{5}$. Namely, instead of the therm $5\omega_{13}^{2}$, there is $5\omega_{15}^{2}.$

\medskip

Also, from \eqref{eq 6} for $x_{2p-1}=0$, $p=3,4,\ldots$ we have
\be\label{eq 8}
\begin{split}
& |\omega_{11}x_{1}+\omega_{31}x_{3}|^{2}+3|\omega_{13}x_{1}+\omega_{33}x_{3}|^{2}
\\
+& 5|\omega_{15}x_{1}+\omega_{35}x_{3}|^{2} + 7|\omega_{17}x_{1}+\omega_{37}x_{3}|^{2}
\leq |x_{1}|^{2}+\frac{|x_{3}|^{2}}{3}.
\end{split}
\ee
From \eqref{eq 8}, for $x_{1}=1$ and  $x_{3}=0 $, since $\omega_{31}=\omega_{13}$, we have
the next inequalities
\[|\omega_{11}|^{2}+3|\omega_{13}|^{2}+5|\omega_{15}|^{2}+7|\omega_{15}|^{2}\leq1, \]
and further
\[
\begin{split}
|\omega_{11}|^{2} &\leq1,\\
|\omega_{11}|^{2}+3|\omega_{13}|^{2}&\leq1, \\
|\omega_{11}|^{2}+3|\omega_{13}|^{2}+5|\omega_{15}|^{2}&\leq1.
\end{split}
\]
This leads to:
\be\label{eq 9}
\begin{split}
|\omega _{11}|& \leq1, \\
|\omega _{13}|&\leq\frac{1}{\sqrt{3}}\sqrt{1-|\omega _{11}|^{2}},\\
|\omega _{15}|&\leq\frac{1}{\sqrt{5}}\sqrt{1-|\omega _{11}|^{2}-3|\omega _{13}|^{2}},\\
|\omega _{17}|&\leq\frac{1}{\sqrt{7}}\sqrt{1-|\omega _{11}|^{2}-3|\omega _{13}|^{2}-5|\omega _{15}|^{2}}.
\end{split}
\ee
We note that we can get the first inequality from \eqref{eq 9} using the fact
\[|a_{2}|=|2\omega_{11}|\leq2 \quad\Rightarrow \quad|\omega_{11}|\leq 1 \] (see \eqref{eq 7}).

\medskip

\section{Proof of Theorem 2}

\medskip

\noindent
\underline{Proof of part $(i)$.} Using the definition of $H_{2}(2)$ given by \eqref{eq 1} and relations \eqref{eq 7}, we have
\[
H_{2}(2)
= 4\omega_{11}\omega_{33}-\frac{7}{3}\omega_{11}^{4}-4\omega_{13}^{2}+4\omega_{11}^{2}\omega_{13}.
\]
NExt, from the fifth relation in \eqref{eq 7} we obtain
\be\label{eq 11}
\omega_{33}=\omega_{15}-\omega_{11}\omega_{13}+\frac{1}{3}\omega_{11}^{3},
\ee
and after combining the two previous relations we have
\[H_{2}(2)=4\omega_{11}\omega_{15}-\omega_{11}^{4}-4\omega_{13}^{2},\]
i.e., 
\[|H_{2}(2)|\leq 4|\omega_{11}||\omega_{15}|+|\omega_{11}|^{4}+4|\omega_{13}|^{2}.\]
Applying  \eqref{eq 9} gives
\be\label{eq 12}
|H_{2}(2)|\leq \frac{4}{\sqrt{5}}|\omega_{11}|\sqrt{1-|\omega _{11}|^{2}-3|\omega _{13}|^{2}} +|\omega_{11}|^{4}+4|\omega_{13}|^{2}:=F_{1}(|\omega _{11}|,|\omega _{13}|),
\ee
where
\be\label{eq 13}
F_{1}(x,y)=\frac{4}{\sqrt{5}}x\sqrt{1-x^{2}-3y^{2}} +x^{4}+4y^{2}.
\ee
Now, we will find the maximum of the function $F_1$ on its domain
\[D_{1}:=\linebreak\left\{0\leq x\leq1,\,0\leq y \leq\frac{1}{\sqrt{3}}\sqrt{1-x^{2}}\right\}.\]

\medskip

Numerically we can verify that the system of equations $\partial F_{1}/\partial x(x,y)=0$ and $\partial F_{1}/\partial y(x,y)=0$ has only one real solution 
$(x_{1},y_{1})=\left(\sqrt{\frac{11}{30}},\frac{1}{30}\sqrt{\frac{281}{2}} \right)=\linebreak (0.60553\ldots ,0.395109\ldots )$ in the interior of $D_{1}$ such that 
then
\[F_{1}(x_{1},y_{1})=1.19889\ldots. \]

Further, let see that maximum values of $F_{1}$ on the boundary of the domain  $D_{1}$.
\begin{itemize}
\item[1)]  For $y=0$, from \eqref{eq 13} we have
\[F_{1}(x,0)=\frac{4}{\sqrt{5}}x\sqrt{1-x^{2}} +x^{4},\quad 0\leq x\leq1. \]
Using the first derivative test we can conclude that the function   $F_{1}(x,0)$  has its maximum at the point
$x_{0}=0.9181\ldots$ which satisfies the equation $5x^8-5x^6+4x^4-4x^2+1=0$ and
\[F_{1}(x_{0},0)=1.3614\ldots. \]
\item[2)]   For $x=0$, since $0\leq y \leq\frac{1}{\sqrt{3}}$, we have
\[F_{1}(0,y)=4y^{2}\leq \frac{4}{3}=1.333\ldots.\]
\item[3)]  Finally,  for $0\leq x\leq1$:
\[F_{1}\left(x,\frac{1}{\sqrt{3}}\sqrt{1-x^{2}}\right)=\frac{1}{3}(3x^{4}-4x^{2}+4)\leq \frac{4}{3}.\]
\end{itemize}

From all the previous facts and \eqref{eq 12} we conclude that $|H_{2}(2)|\leq 1.3614\ldots. $

\bigskip

\noindent
\underline{Proof of part $(ii)$} From the six relation in \eqref{eq 7} and the relation \eqref{eq 11}, after simple calculations, we get
\be\label{eq 14}
\omega_{35}=\omega_{17}-\omega_{11}\omega_{15}+\omega_{11}^{2}\omega_{13}-\omega_{13}^{2}.
\ee
Now, using the relations \eqref{eq 7}, \eqref{eq 11} and \eqref{eq 14}, we obtain
\be\label{eq 15}
\begin{split}
a_{4}&=2(\omega_{15}+3\omega_{11}\omega_{13}+2\omega_{11}^{3}) \\
a_{5}&=2\omega_{17}+6\omega_{11}\omega_{15}+12\omega_{11}^{2}\omega_{13}+3\omega_{13}^{2}+5\omega_{11}^{4}.
\end{split}
\ee
Further, from the definition of $H_{3}(1)$ given by \eqref{eq 2}, the relation \eqref{eq 7} for  $a_{2}$ and $a_{3}$, 
and \eqref{eq 15} for $a_{4}$ and $a_{5}$, after some calculations we have
\[
H_{3}(1)=2\omega_{17}(2\omega_{13}-\omega_{11}^{2})+4\omega_{11}\omega_{13}\omega_{15}
+2\omega_{11}^{3}\omega_{15}-3\omega_{11}^{2}\omega_{13}^{2}-2\omega_{13}^{3} - 4\omega_{15}^2.
\]
So, 
\be\label{eq 16-2}
\begin{split}
|H_{3}(1)| &\leq 2|\omega_{17}||2\omega_{13}-\omega_{11}^{2}|+4|\omega_{11}||\omega_{13}||\omega_{15}|
+2|\omega_{11}|^{3}|\omega_{15}|\\
&+3|\omega_{11}|^{2}|\omega_{13}|^{2}+2|\omega_{13}|^{3} + 4|\omega_{15}|^2.
\end{split}
\ee

\medskip

We satrt analysing the above inequality.

\medskip

Since for the functions from the class $\mathcal{S}$, $ |a_{3}-a_{2}^{2}|\leq1$ (see \cite{duren}), and since from
\eqref{eq 7},
\[ |2\omega_{13}-\omega_{11}^{2}| =|a_{3}-a_{2}^{2}|,\]
we receive
\[
|2\omega_{13}-\omega_{11}^{2}|\leq1.
\]
Using this and the estimate
\[|\omega _{17}|\leq\frac{1}{\sqrt{7}}\sqrt{1-|\omega _{11}|^{2}-3|\omega _{13}|^{2}-5|\omega _{15}|^{2}}
\leq \frac{1}{\sqrt{7}}\sqrt{1-|\omega _{11}|^{2}-3|\omega _{13}|^{2}}\]
given in \eqref{eq 9}, for the first term in \eqref{eq 16-2}, we have
\be\label{eq 18}
2|\omega_{17}||2\omega_{13}-\omega_{11}^{2}|\leq\frac{2}{\sqrt{7}}\sqrt{1-|\omega _{11}|^{2}-3|\omega _{13}|^{2}}.
\ee
Using the estimate for $|\omega _{15}|$ given in \eqref{eq 9} and the estimate in \eqref{eq 18},
inequality  \eqref{eq 16-2} reduces to
\[
\begin{split}
|H_{3}(1)| \leq & \left(\frac{2}{\sqrt{7}}+4|\omega_{11}||\omega_{13}|
+2|\omega_{11}|^{3}\right)\sqrt{1-|\omega _{11}|^{2}-3|\omega _{13}|^{2}}\\
&+\frac{4}{5}-\frac{4}{5}|\omega _{11}|^{2}-\frac{12}{5}|\omega _{13}|^{2}+3|\omega _{11}|^{2}|\omega _{13}|^{2}+2|\omega _{13}|^{3}\\
:= &F_{2}(|\omega _{11}|,|\omega _{13}|),
\end{split}
\]
where
\be\label{eq 19}
F_{2}(x,y)= \left(\frac{2}{\sqrt{7}}+4xy
+2x^{3}\right)\sqrt{1-x^{2}-3y^{2}}
+\frac{4}{5}-\frac{4}{5}x^{2}-\frac{12}{5}y^{2}+3x^{2}y^{2}+2y^{3}
\ee
and $ (x,y)\in D_{1}= \left\{0\leq x\leq1,\,0\leq y \leq\frac{1}{\sqrt{3}}\sqrt{1-x^{2}}\right\}$.

\medskip

Numerical calculation give that the system of equations $\partial F_{2}/\partial x(x,y)=0$ and $\partial F_{2}/\partial y(x,y)=0$ has only two real solutions in the interior of $D_1$, that are  $(x_{2},y_{2})=(0.583\ldots,  0.206\ldots )$ and $(x_{3},y_{3})=(0.0131\ldots, 0.00748\ldots )$ such that
\[F_{2}(x_{2},y_{2})=1.6787\ldots  \quad \text{and}\quad F_{2}(x_{3},y_{3})=1.5559\ldots . \]

\medskip

Now, we consider the maximum values of the function $F_{2}(x,y)$ on the boundary of $D_{1}$.
\begin{itemize}
\item[1)] The relation \eqref{eq 19} for $y=0$ gives 
\[F_{2}(x,0)= \left(\frac{2}{\sqrt{7}}
+2x^{3} \right)\sqrt{1-x^{2}}
+\frac{4}{5}-\frac{4}{5}x^{2},\quad 0\leq x\leq1. \]
Since $ F_{2}(0,0)=\frac{2}{\sqrt{7}}+\frac{4}{5}=1.5559\ldots $,  $ F_{2}(1,0)=0$ and $\partial F_{2}/\partial x(x,0)<0 $ when $0<x<1$, we conclude that
\[F_{2}(x,0)\leq F_{2}(0,0)=\frac{2}{\sqrt{7}}+\frac{4}{5}=1.555928\ldots, \quad 0\leq x\leq1\].

\item[2)] From \eqref{eq 19} for $x=0$ we receive
\[F_{2}(0,y)=\frac{2}{\sqrt{7}}\sqrt{1-3y^{2}}
+\frac{4}{5}-\frac{12}{5}y^{2}+2y^{3},\quad 0\leq y \leq\frac{1}{\sqrt{3}}.\]
Since $F_{2}(0,0)=\frac{2}{\sqrt{7}}+\frac{4}{5}=1.555928\ldots $,\, $F_{2}(0,\frac{1}{\sqrt{3}})=\frac{1}{4}$
and
\[\partial F_{2}/\partial y(0,y) =-\frac{6}{\sqrt{7}}\frac{1}{\sqrt{1-3y^{2}}}-6y(\frac{4}{5}-y)\leq0 \]
when $0\leq y \leq\frac{1}{\sqrt{3}}$, we get 
\[F_{2}(0,y)\leq F_{2}(0,0)=\frac{2}{\sqrt{7}}+\frac{4}{5}=1.5559\ldots,\quad 0\leq y \leq\frac{1}{\sqrt{3}} . \]

\item[3)] At the end,  for $0\leq x\leq1$,
\[F_{2}(x,\frac{1}{\sqrt{3}}\sqrt{1-x^{2}})=\frac{2}{3\sqrt{3}}(1-x^{2})^{\frac{3}{2}}+x^{2}(1-x^{2}).\]
The last function has its maximum $\frac{7}{16}$ for $x=\frac{1}{2}$. So
\[F_{2}\left(x,\frac{1}{\sqrt{3}}\sqrt{1-x^{2}}\right)\leq \frac{7}{16}=0.4375.\]
\end{itemize}

\medskip

Finally, using all the previous facts we conclude that 
\[|H_{3}(1)|\leq 1.6787\ldots . \]

\end{document}